\documentclass[10pt]{amsart}

\usepackage[T1]{fontenc}
\usepackage[utf8]{inputenc}
\usepackage{lmodern}
\usepackage{amsmath,amssymb,mathtools}
\usepackage[protrusion=true,expansion=false]{microtype}
\usepackage{url}
\usepackage[hidelinks]{hyperref}
\allowdisplaybreaks
\theoremstyle{plain}
\newtheorem{theorem}{Theorem}
\newtheorem{lemma}{Lemma}[section]
\newtheorem{proposition}[lemma]{Proposition}

\title[4-gen. non-almost-symmetric nearly Gorenstein semigroup rings]
{Betti sequences of four-generated non-almost-symmetric nearly Gorenstein numerical semigroup rings}

\author{Zhi-Lin Zhang}
\address{Independent Researcher, Taipei, Taiwan}
\email{hsa00000@gmail.com}


\begin{document}

\begin{abstract}
Let \(S\) be a minimally four-generated numerical semigroup that is nearly
Gorenstein but not almost symmetric, and let \(R=\Bbbk[S]\) over an arbitrary
field \(\Bbbk\).  We prove that the defining ideal of \(R\) has four or five
minimal generators when \(R\) has Cohen--Macaulay type two, and exactly six
when it has type three.  Thus the total Betti sequence over the
four-variable polynomial presentation ring is one of
\[
(1,4,5,2),\quad (1,5,6,2),\quad (1,6,8,3),
\]
and all three possibilities occur.  The type-three argument uses a
canonical-ideal exact sequence, a chain-level comparison on Tor, and a
codimension-three Gorenstein reduction.  The resulting classification gives
an affirmative answer, over every field, to Moscariello and Strazzanti's
Question~4.3 on four-generated nearly Gorenstein numerical semigroups that
are not almost symmetric.
\end{abstract}

\maketitle

\section{Introduction}

Let
\[
S=\langle n_1,n_2,n_3,n_4\rangle,
\qquad n_1<n_2<n_3<n_4,
\]
be a minimally generated numerical semigroup, let \(\Bbbk\) be a field, and
write
\[
R=\Bbbk[S]
 =\Bbbk[\mathsf z^{n_1},\ldots,\mathsf z^{n_4}]
 \cong P/I_S,
\qquad
P=\Bbbk[x_1,x_2,x_3,x_4],
\]
where \(\deg x_r=n_r\) and \(x_r\mapsto\mathsf z^{n_r}\).  Put
\(\rho=\mu_P(I_S)\) and \(\tau=\operatorname{type}(S)\).  Herzog's
correspondence identifies \(\rho\) with the cardinality of a minimal
presentation of \(S\) \cite[Proposition~1.5]{Herzog70}, while
Proposition~\ref{prop:betti-identities} gives
\[
(\beta_0^P(R),\beta_1^P(R),\beta_2^P(R),\beta_3^P(R))
=(1,\rho,\rho+\tau-1,\tau).
\]
Thus, in embedding dimension four, the total Betti sequence is determined
once the presentation cardinality \(\rho\) is known.  We study this number
for numerical semigroups that are nearly Gorenstein but not almost
symmetric.

The contrast with the unrestricted four-generated case is sharp.  Bresinsky
showed that the number of minimal defining equations is unbounded among
four-generated numerical semigroups \cite{BresinskyPrime75}.  The nearly
Gorenstein condition, defined through the trace of the canonical module
\cite{HHS19}, imposes much stronger arithmetic constraints: for numerical
semigroups it is equivalent to the existence of an NG-vector
\cite[Definition~1.2 and Proposition~1.1]{MS21}, and in embedding dimension
four it forces \(\tau\le3\) \cite[Theorem~2.4]{MS21}.  Since type one is
symmetric, the non-almost-symmetric case considered here leaves only types
two and three.  This is already a substantial restriction on the final free
module of the resolution, but it does not by itself determine the number of
defining equations.

The exact relation count was left open by Moscariello and Strazzanti.  The
four-generated almost-symmetric case had previously been studied in
\cite{Eto17,HW19}, and Moscariello proved that its type is at most three
\cite[Theorem~1]{Moscariello16}.  In the nearly Gorenstein setting,
Moscariello and Strazzanti obtained the same type bound together with a
normal form for the first non-\(F\) coordinate of an NG-vector, and asked
whether the non-almost-symmetric case always satisfies
\[
\tau=2\Longrightarrow \rho\in\{4,5\},
\qquad
\tau=3\Longrightarrow \rho=6
\]
\cite[Proposition~1.3, Theorem~2.4, and Question~4.3]{MS21}.  More recently,
Chomicz proved the general four-generated estimate
\[
\tau-11\le\rho\le4\tau+5
\]
\cite[Theorem~1.4]{Chomicz26}.  These results constrain the possible size of
a presentation, but they do not determine the minimal relation classes
needed for the exact values above.

What obstructs a direct argument is the passage from semigroup arithmetic
to minimal relations.  An NG-vector gives membership conditions
\[
n_r+f_r-p\in S
\qquad (p\in\operatorname{PF}(S)),
\]
whereas \(\rho\) is the dimension of \(I_S/\mathfrak nI_S\), with
\(\mathfrak n=(x_1,x_2,x_3,x_4)\).  Arithmetic identities readily produce
homogeneous binomials, but they do not decide whether the resulting classes
are minimal or independent.  The same issue persists after replacing
\(I_S\) by an auxiliary ideal: in the type-three argument one obtains
\(I_S=J+(\varphi)\), but this equality alone does not imply
\(\mu_P(I_S)=\mu_P(J)+1\), because adjoining \(\varphi\) can make generators
of \(J\) redundant.  The required comparison must therefore retain
minimality information, not merely exhibit generators.

The proof resolves this mismatch through a translated canonical ideal.  If
\(i\) is the first non-\(F\) coordinate of an NG-vector, the normal form of
\cite[Proposition~1.3]{MS21} selects an index \(l<i\) with
\(f_i=F-n_i+n_l\).  Translating the canonical relative ideal
\[
K(S)=\{a\in\mathbb Z:F-a\notin S\}
\]
by \(n_l\) makes two of its minimal monomial generators equal to the
original variables \(x_l\) and \(x_i\).  The resulting quotient has only two
remaining polynomial variables, which converts the pseudo-Frobenius data
into a low-codimension Gorenstein presentation.  In type two this leads to a
relation-space comparison controlled by the semigroup ideal
\(T_d=\{e\in S:e+d\in S\}\).  In type three it leads instead to a
height-three Gorenstein auxiliary ideal and an exact sequence whose first
Tor map measures precisely whether minimal generators are lost.  These two
comparisons yield the following classification.

\begin{theorem}[Main theorem]\label{thm:main}
Let
\[
S=\langle n_1,n_2,n_3,n_4\rangle,
\quad n_1<n_2<n_3<n_4,
\]
be minimally generated, nearly Gorenstein, and not almost symmetric.  Let
\(\Bbbk\) be any field, let \(R\cong P/I_S\) be the presentation above, and
put \(\rho=\mu_P(I_S)\).  Then
\[
\operatorname{type}(S)=2\Longrightarrow \rho\in\{4,5\},
\quad
\operatorname{type}(S)=3\Longrightarrow \rho=6.
\]
Consequently the total Betti sequence of \(R\) over \(P\) is one of
\[
(1,4,5,2),\quad (1,5,6,2),\quad (1,6,8,3),
\]
and each of these three sequences occurs over every field.
\end{theorem}

The two type cases require different ways of preserving minimality.  For
\(\tau=2\), Section~\ref{sec:type2} identifies the kernel of the relation-space
map arising from the canonical quotient and bounds it through the minimal
\(S\)-generators of \(T_d\), giving \(\rho\in\{4,5\}\).  For \(\tau=3\),
Section~\ref{sec:type3} constructs \(J=I_S\cap Q\) so that \(P/J\) is
Gorenstein and proves, by a chain-level Tor comparison, that the residual
equation does not make any generator of \(J\) redundant.  The
Buchsbaum--Eisenbud structure theorem then reduces the remaining problem to
excluding the three-generator case, which is done by a conormal argument
and a degree-wise binomial graph.  Appendix~\ref{app:sharpness} gives
explicit examples realizing all three sequences.  The argument is valid
over an arbitrary field; it concerns total Betti numbers and uses
essentially the four-generated, non-almost-symmetric setting.

\section{Preliminaries}\label{sec:preliminaries}

Throughout, retain the notation and grading introduced above.  Thus
\(S\) is minimally four-generated, nearly Gorenstein, and not almost
symmetric, \(R=\Bbbk[S]\cong P/I_S\), and \(I_S=\ker\pi\) for
\(\pi(x_r)=\mathsf z^{n_r}\).  We use the same symbol \(x_r\) for its image
in \(R\), and put
\begin{equation*}
F=F(S),\quad \tau=\operatorname{type}(S),\quad
\rho=\mu_P(I_S),\quad
\mathfrak n=(x_1,x_2,x_3,x_4),\quad \mathfrak m=\mathfrak nR.
\end{equation*}
Here \(\mu\) denotes the minimal number of generators, and
\(\mathbb N=\mathbb Z_{\ge0}\).  All graded degrees refer to
\(\deg x_r=n_r\); ordinary total degree will be used only when
explicitly discussing powers of \(\mathfrak n\).
For a finitely generated graded module over a positively graded
\(\Bbbk\)-algebra with degree-zero part \(\Bbbk\), minimal homogeneous
generators are measured by the quotient modulo the homogeneous maximal
ideal.  Localization at that ideal preserves this quotient and minimal
free resolutions.  We use this observation to apply local criteria to
our graded presentations.

Recall that \(F=\max(\mathbb Z\setminus S)\) is the Frobenius number, and
\begin{equation}\label{eq:PF-definition}
\operatorname{PF}(S)=
\{f\in\mathbb Z\setminus S:
  f+s\in S\text{ for all }s\in S\setminus\{0\}\}
\end{equation}
is the set of pseudo-Frobenius numbers of \(S\).  Thus
\(\tau=|\operatorname{PF}(S)|\).  We shall repeatedly use the fact that
\begin{equation}\label{eq:PF-difference}
p-p'\notin S
\quad
\text{for distinct }p,p'\in\operatorname{PF}(S),
\end{equation}
since otherwise \(p=p'+(p-p')\in S\).

Recall that \(R\) is nearly Gorenstein when
\[
\mathfrak m\subseteq\operatorname{tr}_R(\omega_R)
\coloneqq\sum_{\alpha\in\operatorname{Hom}_R(\omega_R,R)}\alpha(\omega_R)
\]
\cite{HHS19}.  We use the equivalent semigroup formulation in terms of
NG-vectors rather than a calculation of this trace ideal.
An \emph{NG-vector} is a tuple
\((f_1,f_2,f_3,f_4)\in\operatorname{PF}(S)^4\) satisfying
\begin{equation}\label{eq:NG-condition}
n_r+f_r-p\in S
\quad
\text{for every }1\le r\le4
\text{ and every }p\in\operatorname{PF}(S).
\end{equation}
By \cite[Definition~1.2 and Proposition~1.1]{MS21}, \(S\) is nearly
Gorenstein if and only if it admits an NG-vector, while \(S\) is almost
symmetric if and only if \((F,F,F,F)\) is an NG-vector.

Since \(S\) is not almost symmetric, every NG-vector has a non-\(F\)
coordinate.  A type-one numerical semigroup is symmetric, hence almost
symmetric; therefore \(\tau\ge2\).  On the other hand,
\cite[Theorem~2.4]{MS21} states that every nearly Gorenstein numerical
semigroup of embedding dimension four has type at most three.  Our
minimal four-generation and nearly Gorenstein assumptions are precisely
its hypotheses.  Thus \(\tau\in\{2,3\}\); the non-almost-symmetric
assumption is what excludes type one and ensures that the first
non-\(F\) coordinate used below exists.  We shall also use the following result of Moscariello and Strazzanti
\cite[Proposition~1.3]{MS21}, stated in our notation.

\begin{lemma}[First non-\(F\) coordinate of an NG-vector]\label{lem:MS-normal}
Every NG-vector has \(f_1=F\).  If \(i\) is its first non-\(F\)
coordinate, then, for some \(l<i\),
\begin{equation*}
f_l=F,\quad f_i=F-d,\quad d=n_i-n_l>0.
\end{equation*}
\end{lemma}

The preceding NG-vector information constrains the pseudo-Frobenius
numbers in each case.  These numbers determine the canonical relative
ideal
\begin{equation}\label{eq:canonical-relative-ideal}
K(S)=\{a\in\mathbb Z:F-a\notin S\}
=\bigcup_{p\in\operatorname{PF}(S)}(F-p+S).
\end{equation}
Its minimal \(S\)-generators are the integers \(F-p\), with
\(p\in\operatorname{PF}(S)\), and the associated monomial fractional
ideal is a graded shift of a canonical module of \(R\).
For the relative-ideal description see \cite[\S1]{MS21}; for the
ring-theoretic identification, \cite[p.~593]{HW19} realizes a canonical
module as \(\sum_{p\in\operatorname{PF}(S)}R\mathsf z^{-p}\).
Multiplication by \(\mathsf z^F\) gives exactly the exponent set
\(K(S)\) above.  Consequently, whenever \(b+K(S)\subseteq S\), the
monomial ideal with exponent set \(b+K(S)\) is a canonical ideal of
\(R\).

We will choose \(b=n_l\): the NG condition at \(l\), where \(f_l=F\),
puts every exponent \(n_l+F-p\) in \(S\), and the exponents corresponding
to \(F\) and \(f_i\) become \(n_l\) and \(n_i\).
To use the resulting quotient, we need both its Gorenstein property and
its socle degree.  These follow from the next elementary description.
Recall that \(a\le_S b\) means \(b-a\in S\).

\begin{lemma}[Socles of monomial quotients]
\label{lem:monomial-socle}
Let \(E\subseteq S\) satisfy \(E+S\subseteq E\), and put
\[
\mathfrak a_E=(\mathsf z^e:e\in E)\subseteq R.
\]
Assume that \(S\setminus E\) is finite and nonempty.  Then the residue classes
\[
\bigl\{\mathsf z^c:
 c\text{ is }\le_S\text{-maximal in }S\setminus E\bigr\}
\]
form a \(\Bbbk\)-basis of
\(\operatorname{Soc}(R/\mathfrak a_E)\).  Consequently,
\(R/\mathfrak a_E\) is Artinian Gorenstein if and only if
\(S\setminus E\) has a unique \(\le_S\)-maximal element.
\end{lemma}

\begin{proof}
The classes of \(\mathsf z^c\), \(c\in S\setminus E\), form a basis.
Multiplication by \(\mathsf z^{n_r}\) sends distinct surviving basis
monomials to distinct basis monomials, so it cannot create cancellations.
Thus a socle element is a linear combination of precisely those basis
monomials for which \(c+n_r\in E\) for every \(r\).  Since the \(n_r\)
generate \(S\) and \(E+S\subseteq E\), this is equivalent to
\(c+s\in E\) for every positive \(s\in S\), or to the stated maximality.
The connected graded Artinian quotient is local, and it is Gorenstein
exactly when its socle has dimension one.
\end{proof}

Two different monomial ideals will occur: the translated canonical ideal
and, in type two, the ideal with exponent set \(T_d\).  Once their
quotients are known to be Gorenstein, the following standard local
criterion gives the same two-variable relation count for either one.

\begin{lemma}[Artinian Gorenstein presentations in at most two variables]
\label{lem:small-presentation}
Let \(D=\Bbbk[y_1,\ldots,y_e]\) be positively graded, with \(0\le e\le2\),
and let \(H\) be a proper homogeneous ideal such that \(D/H\) is
Artinian Gorenstein.  Then \(H\) is generated by a homogeneous regular
sequence of length \(e\), and \(\mu_D(H)=e\).
\end{lemma}

\begin{proof}
Put \(\mathfrak q=(y_1,\ldots,y_e)\).  The local ring
\(D_{\mathfrak q}\) is regular, \(H_{\mathfrak q}\) has height \(e\),
and its quotient is Gorenstein.  For \(e=2\), the codimension-two Gorenstein criterion
\cite[Corollary~21.20]{Eis95} states exactly that such an ideal in a
regular local ring is generated by a regular sequence of length two.  Localization preserves \(H/\mathfrak qH\), so graded
Nakayama gives two homogeneous generators in \(D\).  Their height is
two, hence they form a regular sequence in the Cohen--Macaulay ring
\(D\).  For \(e=1\) the nonzero ideal is principal, and for \(e=0\) it
is zero, generated by the empty regular sequence.
\end{proof}

For a translated canonical ideal, the next proposition also records the
socle degree.  In type three this degree will determine the residual
equation \(\varphi\), and ultimately the complementary degree used in
Gorenstein duality.

\begin{proposition}[Quotients by translated canonical ideals]
\label{prop:quotient-canonical-ideal}
Let \(b\in S\setminus\{0\}\) satisfy \(b+K(S)\subseteq S\), and put
\[
C_b=(\mathsf z^a:a\in b+K(S))\subseteq R.
\]
Then \(C_b\) is a proper canonical ideal, and \(R/C_b\) is Artinian
Gorenstein with socle generated by the class of \(\mathsf z^{F+b}\).
\end{proposition}

\begin{proof}
The canonical-ideal assertion follows from
\eqref{eq:canonical-relative-ideal}.  For \(a\in S\), the definition of
\(K(S)\) gives
\[
a\notin b+K(S)\quad\Longleftrightarrow\quad F+b-a\in S.
\]
Thus every element of \(S\setminus(b+K(S))\) is at most \(F+b\), so
this complement is finite.  Since \(b>0\), one has \(F+b\in S\),
whereas \(F+b\notin b+K(S)\).  The same equivalence shows that every
element of the complement is \(\le_S F+b\).  Hence \(F+b\) is its
unique \(\le_S\)-maximal element.  Lemma~\ref{lem:monomial-socle}
proves the assertions about the quotient and its socle.
\end{proof}

For a finitely generated graded \(P\)-module \(V\), write
\[
\beta_{p,e}^P(V)
=\dim_{\Bbbk}\operatorname{Tor}_p^P(V,\Bbbk)_e,
\quad
\beta_p^P(V)=\sum_e\beta_{p,e}^P(V).
\]
We also record the Betti-number identity used in both type cases.

\begin{proposition}\label{prop:betti-identities}
The minimal \(P\)-resolution of \(R\) has length three, and
\[
(\beta_0^P(R),\beta_1^P(R),\beta_2^P(R),\beta_3^P(R))
=(1,\rho,\rho+\tau-1,\tau).
\]
\end{proposition}

\begin{proof}
Let \(F_\bullet\to R\) be the minimal graded \(P\)-free resolution.
The ring \(R\) is a one-dimensional domain, so \(R_{\mathfrak m}\) is
Cohen--Macaulay.  It is finite over the regular local ring
\(P_{\mathfrak n}\) of dimension four and has finite projective dimension.
The Auslander--Buchsbaum formula therefore gives
\(\operatorname{pd}_P R=3\), using preservation of minimality under
localization \cite[Theorem~19.9]{Eis95}.
By \cite[Corollary~21.16]{Eis95}, dualizing the localized minimal
resolution gives a minimal resolution of the canonical module.  Hence
\[
\beta_3^P(R)=\operatorname{rank}F_3
=\mu(\omega_{R_{\mathfrak m}})=\tau.
\]
Here \(\mu(\omega_{R_{\mathfrak m}})\) is the Cohen--Macaulay type by
\cite[Exercise~21.14(a)]{Eis95}, and for a numerical semigroup ring this
type equals \(|\operatorname{PF}(S)|=\tau\) \cite[\S1]{MS21}.
The first two free ranks are \(1,\rho\).  Finally, \(I_S\ne0\) implies
\(R\otimes_P\operatorname{Frac}(P)=0\); the alternating sum of free
ranks is therefore zero.  Thus
\(1-\rho+\beta_2^P(R)-\tau=0\), as required.
\end{proof}

\section{The type-two case}\label{sec:type2}

The quotient by \(C=(x_l,x_i)R\) is useful because it has only two
variables and hence exactly two minimal equations.  We must bound the
relation classes that disappear in this quotient, not count its equations
alone.  The two lemmas below separate this task into a relation-space
comparison and a bound on minimal semigroup-ideal generators.

Assume throughout this section that \(\tau=2\).  Choose an NG-vector and
let \(i\) be its first non-\(F\) coordinate.  Lemma~\ref{lem:MS-normal}
supplies an index \(l<i\) such that
\(\operatorname{PF}(S)=\{F,F-d\}\), where \(d=n_i-n_l>0\).

Let \(\{j,k\}=\{1,2,3,4\}\setminus\{i,l\}\), with their order
immaterial, and set \(T_d=\{e\in S:e+d\in S\}\).
Let \(\operatorname{Min}_S(T_d)\) denote the set of elements minimal for
the order \(a\le_S b\) if \(b-a\in S\).  Since
\(F-d\in\operatorname{PF}(S)\), one has \(d\notin S\); otherwise
\(F=(F-d)+d\in S\).  In particular, \(0\notin T_d\).

The set \(T_d\) is dictated by equality of degrees in a binomial
\(x_lU-x_iV\): if \(\deg V=e\), then \(\deg U=e+d\).
The next lemma turns this observation into a bound on independent
relation classes.  Its subtraction of one is essential for the final
bound \(\rho\le5\).

\begin{lemma}[Comparison with \(T_d\)]\label{lem:type2-comparison}
One has
\[
\rho-2\le |\operatorname{Min}_S(T_d)|-1.
\]
\end{lemma}

\begin{proof}
Since \(K(S)=S\cup(d+S)\) and \(n_i=n_l+d\), the ideal
\(C=(\mathsf z^{n_l},\mathsf z^{n_i})R=(x_l,x_i)R\) has exponent set
\(n_l+K(S)\).  Thus it is the canonical ideal \(C_{n_l}\) of
Proposition~\ref{prop:quotient-canonical-ideal}.
Put \(Q=(x_l,x_i)\), \(\bar P=P/Q\),
\(\bar I=(I_S+Q)/Q\), and
\(\bar{\mathfrak n}=\mathfrak n\bar P\).
Proposition~\ref{prop:quotient-canonical-ideal} makes
\(R/C\cong\bar P/\bar I\) Artinian Gorenstein, and
Lemma~\ref{lem:small-presentation} gives \(\mu_{\bar P}(\bar I)=2\), since
\(\bar P\cong\Bbbk[x_j,x_k]\).
Hence the natural surjection
\(I_S/\mathfrak nI_S\to
\bar I/\bar{\mathfrak n}\bar I\)
has kernel
\(((I_S\cap Q)+\mathfrak nI_S)/\mathfrak nI_S\)
of dimension \(\rho-2\).  In particular,
\begin{equation}\label{eq:type2-relation-kernel}
\rho-2=
\dim_{\Bbbk}\frac{(I_S\cap Q)+\mathfrak nI_S}{\mathfrak nI_S}.
\end{equation}

For \(e\in T_d\), choose monomials \(U_e,V_e\) with
\(\deg U_e=e+d\) and \(\deg V_e=e\), and set
\(B_e=x_lU_e-x_iV_e\).
Since \(n_l+e+d=n_i+e\), one has \(B_e\in I_S\cap Q\).
The classes of these binomials span the kernel.  Indeed, each
homogeneous component of \(I_S\cap Q\) is spanned by binomials of
equal \(S\)-degree.  If both terms are divisible by \(x_l\), or both
by \(x_i\), then the binomial lies in \(\mathfrak nI_S\).
Otherwise, after interchanging the terms, it has the form
\(x_lU-x_iV\); if \(e=\deg V\), equality of degrees gives
\(\deg U=e+d\), so \(e\in T_d\).

The class of \(B_e\) depends only on \(e\), since changing \(U_e\) or
\(V_e\) changes \(B_e\) by an element of
\(QI_S\subseteq\mathfrak nI_S\).
Moreover, if \(e=e'+s\) with \(e'\in T_d\) and
\(s\in S\setminus\{0\}\), choose a monomial \(W\) of degree \(s\)
and take \(U_e=WU_{e'}\), \(V_e=WV_{e'}\).
Then \(B_e=WB_{e'}\in\mathfrak nI_S\).
Thus only elements of \(\operatorname{Min}_S(T_d)\) can contribute
nonzero classes.

Finally, \(n_l\in T_d\) because \(n_l+d=n_i\), and it is minimal in
\(T_d\): otherwise \(n_l=e+s\), with \(e\in T_d\) and
\(s\in S\setminus\{0\}\), would express the minimal generator \(n_l\)
as a sum of two positive elements of \(S\).
For \(e=n_l\), however, we may take
\(U_{n_l}=x_i\) and \(V_{n_l}=x_l\), so \(B_{n_l}=0\).
Consequently
\[
\rho-2\le
\bigl|\operatorname{Min}_S(T_d)\setminus\{n_l\}\bigr|
=|\operatorname{Min}_S(T_d)|-1.
\]
\end{proof}

It remains to bound the minimal generators of \(T_d\).
The useful change of viewpoint is to regard it as the exponent set of a
monomial ideal: its quotient has a one-dimensional socle.  If this ideal
contains at least two of the original variables, the two-variable
Gorenstein criterion applies.  If it contains only one, that criterion
no longer suffices; the NG-vector covering then forces all but at most
one minimal element into an explicit set of three elements.  This
explains the two cases in the next proof.

\begin{lemma}[Bounding the minimal elements of \(T_d\)]\label{lem:type2-Td-bound}
One has
\[
|\operatorname{Min}_S(T_d)|\le4.
\]
\end{lemma}

\begin{proof}
Put \(\mathcal T=T_d\cup(d+T_d)\).  Since \(T_d+S=T_d\), one has
\(\mathcal T+S=\mathcal T\).  For every minimal generator \(n_r\),
the coordinate \(f_r\) is either \(F\) or \(F-d\).  If \(f_r=F\),
\eqref{eq:NG-condition} with \(p=F-d\) gives \(n_r+d\in S\), so
\(n_r\in T_d\).  If \(f_r=F-d\), the same condition with \(p=F\)
gives \(n_r-d\in S\).  Since \((n_r-d)+d=n_r\in S\), it follows
that \(n_r-d\in T_d\), so \(n_r\in d+T_d\).
Thus \(\mathcal T\) contains every minimal generator, and its stability
under addition by \(S\) gives \(S\setminus\{0\}\subseteq\mathcal T\).
The reverse inclusion follows from \(d>0\), \(0\notin T_d\), and the
definition of \(T_d\).  Consequently
\begin{equation}\label{eq:type2-Td-union}
S\setminus\{0\}=T_d\cup(d+T_d).
\end{equation}
Set \(\mathfrak a_{T_d}=(\mathsf z^e:e\in T_d)\subseteq R\).

We first show that \(R/\mathfrak a_{T_d}\) is Artinian Gorenstein.
The set \(S\setminus T_d\) is finite, since every sufficiently large
integer belongs to \(T_d\), and it is nonempty because \(0\notin T_d\).
Let \(c\) be \(\le_S\)-maximal in \(S\setminus T_d\).
Then \(c+d\notin S\), while maximality gives \(c+n_r\in T_d\) for
every \(r\).  Hence \(c+d+n_r\in S\) for every minimal generator
\(n_r\), and therefore \(c+d+s\in S\) for every
\(s\in S\setminus\{0\}\).  Thus
\(c+d\in\operatorname{PF}(S)=\{F,F-d\}\).
The value \(F\) is impossible, since it would give
\(c=F-d\notin S\).  Hence every \(\le_S\)-maximal element of
\(S\setminus T_d\) is the same integer \(F-2d\).  Since the finite
nonempty set \(S\setminus T_d\) has at least one maximal element, it
has exactly one.  Lemma~\ref{lem:monomial-socle}, applied to
\(E=T_d\), now shows that \(R/\mathfrak a_{T_d}\) is Artinian
Gorenstein.

Put \(U=\{r:n_r\in T_d\}\).  Since \(f_1=F\), the NG inequality
against \(F-d\) gives \(n_1+d\in S\), so \(1\in U\).
Suppose first that \(|U|\ge2\).  Let
\(\widetilde{\mathfrak a}=\pi^{-1}(\mathfrak a_{T_d})\subseteq P\).
For \(r\in U\), the image of \(x_r\) lies in \(\mathfrak a_{T_d}\), so
\(x_r\in\widetilde{\mathfrak a}\).  Hence
\(R/\mathfrak a_{T_d}\cong P/\widetilde{\mathfrak a}\) is presented by
the polynomial ring in the variables \(x_r\) with \(r\notin U\), and
its defining ideal is generated by \(4-|U|\) elements by
Lemma~\ref{lem:small-presentation}.  Here the quotient has just been
proved Artinian Gorenstein, and \(4-|U|\le2\), so all hypotheses of
that lemma hold.  Lifting these generators to \(P\), together with the
variables \(x_r\) for \(r\in U\), generates
\(\widetilde{\mathfrak a}\) with at most
\(|U|+(4-|U|)=4\) elements.  Passing to \(R\), we obtain
\(|\operatorname{Min}_S(T_d)|
=\mu_R(\mathfrak a_{T_d})\le4\).

It remains to consider \(U=\{1\}\).  For \(r=2,3,4\), one has
\(n_r\notin T_d\), so \eqref{eq:type2-Td-union} gives
\(n_r-d\in T_d\).  Moreover \(f_r\ne F\): otherwise the NG inequality
against \(F-d\) would give \(n_r+d\in S\), hence \(n_r\in T_d\), a
contradiction.  Thus the first non-\(F\) coordinate is \(2\), and
Lemma~\ref{lem:MS-normal} gives \(d=n_2-n_1\).  Consequently
\(D=\{n_1,n_2-d,n_3-d,n_4-d\}
=\{n_1,n_3-d,n_4-d\}\) has at most three elements.

Let \(e\in\operatorname{Min}_S(T_d)\setminus D\).
If a factorization of \(e+d\) used some \(n_r\) with \(r\ge2\), say
\(e+d=n_r+s\), then \(e=(n_r-d)+s\).
If \(s=0\), this gives \(e=n_r-d\in D\); if \(s\ne0\), then
\(n_r-d\in T_d\), contradicting the \(S\)-minimality of \(e\).
Hence \(e+d=a n_1\) for some positive integer \(a\).
If also \(e'+d=a'n_1\) with \(a<a'\), then
\(e'-e=(a'-a)n_1\in S\), so \(e\le_S e'\).
Thus at most one element of
\(\operatorname{Min}_S(T_d)\setminus D\) can occur, and therefore
\(|\operatorname{Min}_S(T_d)|\le |D|+1\le4\).
\end{proof}

\begin{proposition}\label{prop:type2}
If \(\tau=2\), then \(\rho\in\{4,5\}\).
\end{proposition}

\begin{proof}
Since \(\operatorname{ht} I_S=3\), one has \(\rho\ge3\).
If \(\rho=3\), then \(I_S\) is a complete intersection, so \(R\) is
Gorenstein and has type one, contrary to \(\tau=2\).  Hence \(\rho\ge4\).

Lemmas~\ref{lem:type2-comparison} and~\ref{lem:type2-Td-bound} give
\[
\rho-2\le |\operatorname{Min}_S(T_d)|-1\le3,
\]
so \(\rho\le5\).
\end{proof}

\section{The type-three case}\label{sec:type3}

Assume throughout this section that \(\tau=3\).  Unlike the type-two
argument, the proof will not count the kernel of a two-variable
presentation directly.  Instead it constructs a height-three Gorenstein
ideal \(J\) and proves that adjoining one equation gives \(I_S\) without
making any generator of \(J\) redundant.  The three tasks are therefore
\[
\text{construct \(B=P/J\) Gorenstein}
\ \longrightarrow\ \rho=\mu_P(J)+1
\ \longrightarrow\ \mu_P(J)=5.
\]
The first transition requires a Tor-map calculation; the second requires
both the Pfaffian theorem and a separate exclusion of the
complete-intersection alternative.  Neither transition follows merely
from the Gorenstein property of \(B\).

Choose an NG-vector and let \(i\) be its first non-\(F\) coordinate.
Lemma~\ref{lem:MS-normal} supplies an index \(l<i\).  Set
\begin{equation}\label{eq:type3-PF-data}
f_i=F-d,\quad d=n_i-n_l>0,\quad
\operatorname{PF}(S)=\{F,f_i,g\},
\end{equation}
where \(f_l=F\) and \(g\) is the remaining pseudo-Frobenius number.
Let \(\{j,k\}=\{1,2,3,4\}\setminus\{l,i\}\).

\subsection{Lifting the canonical ideal and constructing the fibre product}

The translate \(n_l+K(S)\) contains the generators \(n_l,n_i\) and one
further exponent \(n_l+F-g\).  We need a representative \(M\) of the
third exponent involving neither \(x_l\) nor \(x_i\), so that
\(x_l,x_i,M\) is a regular sequence in \(P\).  We also choose its
orientation so that \(x_j\mid M\) and \(g+n_j=\lambda n_k\); this
will make a pure power of \(x_k\) survive modulo \(Q\) and supply an
explicit residual equation.  The normal form
\(g=\lambda n_k-n_j\) is part of \cite[Theorem~2.4]{MS21} under the
present nearly Gorenstein, non-almost-symmetric, four-generated
hypotheses.  The following argument recalls the relevant arithmetic
while making this compatible choice of \(M\) explicit.

\begin{proposition}[Canonical ideal construction in type three]\label{prop:type3-canonical}
After possibly interchanging \(j\) and \(k\), there exist integers
\(\theta,\lambda\ge1\) and \(\nu\ge0\) such that, with
\(M=x_j^\theta x_k^\nu\),
\begin{equation}\label{eq:type3-monomial}
\deg M=n_l+F-g,\quad g+n_j=\lambda n_k.
\end{equation}
The exponent set of
\(C=(x_l,x_i,M)R\) is \(n_l+K(S)\).  In particular,
\(C\) is a canonical ideal minimally generated by \(x_l,x_i,M\), and
\(A_C=R/C\) is Artinian Gorenstein.  Moreover, for every \(e\in S\),
\begin{equation}\label{eq:type3-C-membership}
\mathsf z^e\in C
\quad\Longleftrightarrow\quad
F+n_l-e\notin S.
\end{equation}
\end{proposition}

\begin{proof}
By \eqref{eq:NG-condition}, with \(r=l\) and \(p=g\), and since
\(f_l=F\), one has \(n_l+F-g\in S\).  Choose
\(a_l,a_i,a_j,a_k\in\mathbb Z_{\ge0}\) such that
\[
n_l+F-g=a_ln_l+a_in_i+a_jn_j+a_kn_k.
\]
If \(a_l>0\), then \(F-g\in S\), contrary to
\eqref{eq:PF-difference}.  If \(a_i>0\), then by
\eqref{eq:type3-PF-data},
\[
n_l+F-g-n_i=F-g-d=f_i-g\in S,
\]
again contradicting \eqref{eq:PF-difference}.  Hence
\(a_l=a_i=0\).  Since \(n_l+F-g>0\), at least one of \(a_j,a_k\)
is positive.  After possibly interchanging \(j\) and \(k\), write
\(n_l+F-g=\theta n_j+\nu n_k\), where
\(\theta\ge1\) and \(\nu\ge0\), and set
\(M=x_j^\theta x_k^\nu\).

Since \(g\in\operatorname{PF}(S)\), \eqref{eq:PF-definition} gives
\(g+n_j\in S\).  In any factorization of \(g+n_j\), the coefficient
of \(n_j\) is zero, since otherwise \(g\in S\).  The coefficients of
\(n_l\) and \(n_i\) are also zero: otherwise, together with
\((\theta-1)n_j+\nu n_k\in S\), one would obtain respectively
\[
F=(g+n_j-n_l)+(\theta-1)n_j+\nu n_k\in S
\]
or
\[
f_i=(g+n_j-n_i)+(\theta-1)n_j+\nu n_k\in S,
\]
both impossible.  Hence \(g+n_j=\lambda n_k\) for some
\(\lambda\ge0\).  The case \(\lambda=0\) would give \(g=-n_j\),
and then \(g+n_k=n_k-n_j\in S\), contrary to minimal generation.
Thus \(\lambda\ge1\).

We now identify the exponent set of \(C\). By \eqref{eq:canonical-relative-ideal} and
\eqref{eq:type3-PF-data},
\[
K(S)=S\cup(d+S)\cup(F-g+S).
\]
Since \(n_l+d=n_i\) and \(\deg M=n_l+F-g\),
\[
n_l+K(S)
=(n_l+S)\cup(n_i+S)\cup(\deg M+S),
\]
which is precisely the exponent set of \(C=(x_l,x_i,M)R\).
Thus \(C\) is canonical.  Its generators are minimal because their
exponents are \(n_l+F-p\), with \(p\in\{F,f_i,g\}\), and the
difference of two such exponents is the difference of two distinct
pseudo-Frobenius numbers, which is not in \(S\) by
\eqref{eq:PF-difference}.  Applying
Proposition~\ref{prop:quotient-canonical-ideal} with \(b=n_l\) shows
that \(A_C=R/C\) is Artinian Gorenstein.  Finally, for \(e\in S\),
\[
\mathsf z^e\in C
\quad\Longleftrightarrow\quad
e-n_l\in K(S)
\quad\Longleftrightarrow\quad
F+n_l-e\notin S,
\]
which proves \eqref{eq:type3-C-membership}.
\end{proof}

Fix the monomial \(M\) constructed above and set
\[
Q=(x_l,x_i,M)\subseteq P,\quad
A_Q=P/Q,\quad
J=I_S\cap Q,\quad
B=P/J.
\]
Thus \(J\) is the auxiliary ideal obtained by intersecting the defining
ideal \(I_S\) with the complete-intersection ideal \(Q\); we shall
compare its minimal generators with those of \(I_S\).  Since \(QR=C\),
\[
A_C=R/C\cong P/(I_S+Q).
\]
Moreover, \(x_l,x_i,M\) is a regular sequence: after killing the two
variables, \(M\) is a nonzero monomial in the domain
\(\Bbbk[x_j,x_k]\).  Thus \(A_Q\) is a one-dimensional complete
intersection.  For Ogoma's theorem, the second quotient map
\(A_Q\to A_C\) must also have a canonical kernel.  It is enough that
this kernel be principal and generated by a nonzerodivisor, since
\(A_Q\) is Gorenstein.  The next lemma proves exactly this, using the
two-variable presentation of \(A_C\).  The socle degree fixes the
degree \(q\) of the residual equation; retaining that degree is
essential for the later duality calculation.

\begin{lemma}[The comparison exact sequence]\label{lem:type3-equation}
Set \(q=g+n_j+n_k=(\lambda+1)n_k\).  Choose a monomial \(L\in P\)
of degree \(g+n_k\), and put
\(\varphi=x_k^{\lambda+1}-x_jL\).  Then \(\varphi\in I_S\), and its
image \(\bar\varphi\) in \(A_Q=P/Q\) is a homogeneous
nonzerodivisor of degree \(q\) satisfying
\[
(I_S+Q)/Q=(\bar\varphi).
\]
Consequently,
\[
I_S=J+(\varphi),
\quad
J:\varphi=Q,
\]
and multiplication by \(\varphi\) gives an exact sequence
\begin{equation}\label{eq:type3-exact-sequence}
0\longrightarrow A_Q(-q)\xrightarrow{\ \cdot\varphi\ }B
\longrightarrow R\longrightarrow0.
\end{equation}
\end{lemma}

\begin{proof}
To use the two-variable presentation of \(A_C\), let
\(D=P/(x_l,x_i)\cong\Bbbk[x_j,x_k]\) and
\(H=(I_S+Q)/(x_l,x_i)\).  The quotient \(A_C\cong D/H\) is
Artinian Gorenstein by Proposition~\ref{prop:type3-canonical}, and
\(D\) is a positively graded polynomial ring in two variables.
Lemma~\ref{lem:small-presentation} therefore gives two minimal
homogeneous generators of \(H\).  The class of \(M\) is one of
them.  Indeed, if \(M\in(x_j,x_k)H\), then lifting to \(P\) and
passing to \(R\) gives
\(M\in(x_l,x_i)R+\mathfrak m C\).  Since
\(C=(x_l,x_i,M)R\), graded Nakayama's lemma would imply
\(C=(x_l,x_i)R\), contradicting the minimality of the generators of
\(C\).  Hence \(H=(M,\bar\psi)\) for some homogeneous
\(\psi\in I_S\).  As \(H\) has height two, \(M,\bar\psi\) is a
regular sequence in \(D\).  Equivalently, the image of \(\psi\) in
\(A_Q\) is a nonzerodivisor and generates \((I_S+Q)/Q\).

The degree of the second generator is forced by the socle of \(A_C\).
By Proposition~\ref{prop:quotient-canonical-ideal} with \(b=n_l\), this
socle is generated by the class of \(\mathsf z^{F+n_l}\).
As \(D=\Bbbk[x_j,x_k]\) has weights \(n_j,n_k\) and
\(M,\bar\psi\) is a homogeneous regular sequence, its Hilbert series is
\[
H_{A_C}(t)
=
\frac{(1-t^{\deg M})(1-t^{\deg\psi})}
     {(1-t^{n_j})(1-t^{n_k})}.
\]
Since \(A_C\) is Artinian, this is a polynomial whose highest nonzero
degree is \(\deg M+\deg\psi-n_j-n_k\).  Because \(A_C\) is
Gorenstein, its highest nonzero graded component is its socle; thus its
socle degree is
\(\deg M+\deg\psi-n_j-n_k\).  Equivalently, the weighted complete-intersection formula
\(\omega_{D/H}\cong(D/H)(\deg M+\deg\psi-n_j-n_k)\)
applies here: the variables have the positive weights \(n_j,n_k\), and
\(M,\bar\psi\) is a homogeneous regular sequence
\cite[Exercise~21.16]{Eis95}.  Hence
\[
F+n_l=\deg M+\deg\psi-n_j-n_k.
\]
Using \eqref{eq:type3-monomial}, this gives
\[
\deg\psi=g+n_j+n_k=q.
\]
Since \(g\in\operatorname{PF}(S)\), one has \(g+n_k\in S\), so the
monomial \(L\) exists.  Moreover, by \eqref{eq:type3-monomial},
\[
\deg(x_k^{\lambda+1})
=(\lambda+1)n_k
=q
=\deg(x_jL).
\]
Both monomials map to \(\mathsf z^q\) in \(R\), hence
\(\varphi\in I_S\).  Its image in \(A_Q\) is nonzero because \(Q\)
is monomial and \(x_k^{\lambda+1}\notin Q\).  Since
\((I_S+Q)/Q=(\bar\psi)\), the ring \(A_Q\) is positively graded with
\((A_Q)_0=\Bbbk\), and \(\bar\varphi,\bar\psi\) have the same degree,
we have \(\bar\varphi=c\bar\psi\) for some
\(c\in\Bbbk^\times\).  Thus \(\bar\varphi\) is a nonzerodivisor and
generates \((I_S+Q)/Q\).

For \(f\in I_S\), the equality
\((I_S+Q)/Q=(\bar\varphi)\) gives
\(f-a\varphi\in Q\) for some \(a\in P\).  Since
\(f,\varphi\in I_S\), this difference belongs to
\(I_S\cap Q=J\), and therefore \(I_S=J+(\varphi)\).
Since \(\bar\varphi\) is a nonzerodivisor in \(A_Q=P/Q\),
\(Q:\varphi=Q\), while \(\varphi\in I_S\) gives
\(I_S:\varphi=P\).  Hence
\[
J:\varphi
=(I_S\cap Q):\varphi
=(I_S:\varphi)\cap(Q:\varphi)
=Q.
\]
Finally, \(Q\varphi\subseteq J\), so multiplication by \(\varphi\)
induces \(A_Q(-q)\to B\).  Its kernel is
\((J:\varphi)/Q=0\), and its image is
\((J+(\varphi))/J=I_S/J\).  The quotient is \(P/I_S=R\), proving
\eqref{eq:type3-exact-sequence}.
\end{proof}

The exact sequence alone does not make \(B\) Gorenstein.  The reason for
intersecting with \(Q\), rather than choosing an arbitrary ideal contained
in \(I_S\), is that it realizes \(B\) as a fibre product with the two
canonical kernels just established.

\begin{proposition}[The auxiliary Gorenstein ring]\label{prop:type3-gorenstein}
The ring \(B=P/J\) is a one-dimensional graded Gorenstein ring.
\end{proposition}

\begin{proof}
Consider the two quotient maps \(R\to A_C\) and \(A_Q\to A_C\).
The map \(p\mapsto(p+I_S,p+Q)\) identifies
\(B=P/(I_S\cap Q)\) with \(R\times_{A_C}A_Q\): its kernel is the
intersection, and two cosets with the same image modulo \(I_S+Q\)
have a common lift in \(P\).  Since \(A_C\) is Artinian local,
localization at the homogeneous maximal ideals gives
\[
B_{\mathfrak nB}
\cong R_{\mathfrak m}\times_{A_C}(A_Q)_{\mathfrak nA_Q}.
\]

We use the sufficient direction of Ogoma's theorem
\cite[Theorem~4]{Ogoma87}.  For surjections of Noetherian local rings
\(A_1\to A_0\leftarrow A_2\), with nonzero kernels
\(\mathfrak a_1,\mathfrak a_2\), suppose that \(A_1\) is
Cohen--Macaulay and \(\operatorname{grade}_{A_1}\mathfrak a_1\ge1\),
and that \(A_2\) is equidimensional and satisfies \((S_1)\).
The theorem makes the fibre product Gorenstein if \(A_2\) is also
Cohen--Macaulay and each \(\mathfrak a_r\) is a canonical module of
\(A_r\).  Only this implication, not the converse classification, is
needed here.

Here \(A_1=R_{\mathfrak m}\),
\(A_2=(A_Q)_{\mathfrak nA_Q}\), and \(A_0=A_C\).  Both quotient maps
are surjective, and their kernels are nonzero, so this is the
nontrivial case.  The ring \(A_1\) is a one-dimensional
Cohen--Macaulay local ring.  Its kernel \(C_{\mathfrak m}\) is a
canonical module by Proposition~\ref{prop:type3-canonical}; since it
contains the nonzerodivisor \(x_l\), it has grade at least one.  The
second kernel is
\((\bar\varphi)(A_Q)_{\mathfrak nA_Q}\).  By
Lemma~\ref{lem:type3-equation}, \(\bar\varphi\) is a nonzerodivisor,
so multiplication by \(\bar\varphi\) identifies this kernel with
\((A_Q)_{\mathfrak nA_Q}\).  The latter is a one-dimensional local
complete intersection, hence Gorenstein, Cohen--Macaulay,
equidimensional, and satisfies \((S_1)\); it is therefore a canonical
module over itself.  All of Ogoma's hypotheses are thus satisfied, and
\(B_{\mathfrak nB}\) is Gorenstein.

The support of \(P/(I_S\cap Q)\) is the union of those of \(P/I_S\)
and \(P/Q\), so \(\dim B=1\).  The depth lemma applied to
\eqref{eq:type3-exact-sequence} makes \(B\) Cohen--Macaulay.
Its graded canonical module \(\omega_B\) has one minimal generator,
since its localization at \(\mathfrak nB\) does.  Choose the resulting
homogeneous surjection \(B(-a)\twoheadrightarrow\omega_B\).
Locally this is an isomorphism, as a surjection from a local ring to a
module isomorphic to that ring.  Its finitely generated graded kernel
therefore vanishes after localization at \(\mathfrak nB\), and hence
vanishes by graded Nakayama.  Thus \(\omega_B\cong B(-a)\), proving
that \(B\) is graded Gorenstein.
\end{proof}

\subsection{The Tor comparison and preservation of minimal generators}
\label{subsec:tor-comparison}

The equality \(I_S=J+(\varphi)\) gives a generating set but does not
preserve minimality automatically.  Apply \(\operatorname{Tor}^P(-,\Bbbk)\)
to \eqref{eq:type3-exact-sequence}, and write
\[
\alpha_p:\operatorname{Tor}_p^P(A_Q(-q),\Bbbk)
\longrightarrow\operatorname{Tor}_p^P(B,\Bbbk)
\]
for the maps induced by multiplication by \(\varphi\).
Our target is \(\alpha_1=0\).  Two of its three source classes admit
direct binomial decompositions; the third needs a complementary-degree
vanishing.  We first extract the high-degree Tor information needed for
that vanishing, before returning to the three classes themselves.

\begin{lemma}[Higher Tor and the duality degree]
\label{lem:type3-higher-tor}
Put \(s_P=n_1+n_2+n_3+n_4\) and
\[
\ell_B=q+n_l+n_i+\deg M=F+n_l+s_P.
\]
The last free module of the minimal resolution of \(B\) is
\(P(-\ell_B)\).  Moreover \(\alpha_2=0\), there is a graded injection
\begin{equation}\label{eq:type3-tor2-injection}
\operatorname{Tor}_2^P(B,\Bbbk)
\hookrightarrow\operatorname{Tor}_2^P(R,\Bbbk),
\end{equation}
and, for every integer \(e\),
\begin{equation}\label{eq:type3-duality}
\beta_{1,e}^P(B)=\beta_{2,\ell_B-e}^P(B).
\end{equation}
\end{lemma}

\begin{proof}
Both \(R\) and \(B\) are one-dimensional Cohen--Macaulay quotients of
\(P\), so their projective dimensions are three, by the same
Auslander--Buchsbaum calculation as in
Proposition~\ref{prop:betti-identities}.  The top free rank for \(B\) is
one because \(B\) is Gorenstein.  The Koszul resolution of
\(A_Q(-q)\), on the regular sequence \(x_l,x_i,M\), has top degree
\(\ell_B\) and free ranks \((1,3,3,1)\).  Thus the upper part of the
long exact sequence is
\[
\begin{aligned}
0&\longrightarrow\Bbbk(-\ell_B)
 \xrightarrow{\alpha_3}\operatorname{Tor}_3^P(B,\Bbbk)
 \longrightarrow\operatorname{Tor}_3^P(R,\Bbbk)\\
 &\xrightarrow{\delta_3}\operatorname{Tor}_2^P(A_Q(-q),\Bbbk)
 \xrightarrow{\alpha_2}\operatorname{Tor}_2^P(B,\Bbbk)
 \longrightarrow\operatorname{Tor}_2^P(R,\Bbbk).
\end{aligned}
\]
The injection \(\alpha_3\) is an isomorphism between one-dimensional
spaces, fixing the top degree of \(B\).  Exactness then makes
\(\delta_3\) injective.  Its source has dimension \(\tau=3\), and its
target has dimension three by the Koszul resolution, so it is an
isomorphism.  Therefore \(\alpha_2=0\), which gives
\eqref{eq:type3-tor2-injection}.  This dimension match is the precise
use of the type-three hypothesis in the comparison.

Write \(G_\bullet\) for the minimal graded resolution of \(B\).
Duality applies to the Cohen--Macaulay quotient \(B\) of the regular
ring \(P\), with \(\omega_P=P(-s_P)\).
Dualizing \(G_\bullet\) into \(P(-s_P)\) resolves
\(\omega_B\) \cite[Corollary~21.16 and \S21.11]{Eis95}.
The zeroth free module of the dual resolution is
\(P(\ell_B-s_P)\).  Since \(B\) is graded Gorenstein, this identifies
\(\omega_B\cong B(\ell_B-s_P)\).  Shifting back and using uniqueness
of minimal graded resolutions gives
\(G_2\cong\operatorname{Hom}_P(G_1,P(-\ell_B))\), which is
\eqref{eq:type3-duality}.
\end{proof}

The third class has degree \(q+\deg M\), whose complement under
\eqref{eq:type3-duality} is
\(\ell_B-(q+\deg M)=n_l+n_i\).  This is a degree where minimal
generation of the semigroup makes the required homology especially
simple.  The vanishing for \(R\) in the next lemma uses only minimal
four-generation, not the nearly Gorenstein condition.

\begin{lemma}[Vanishing in the complementary degree]
\label{lem:type3-complementary-vanishing}
One has \(\beta_{2,n_l+n_i}^P(R)=0\), and consequently
\(\beta_{1,q+\deg M}^P(B)=0\).
\end{lemma}

\begin{proof}
For \(e\in S\), set
\[
\Delta_e=\left\{A\subseteq\{1,2,3,4\}:
 e-\sum_{r\in A}n_r\in S\right\}.
\]
The squarefree-divisor-complex formula of Campillo and Mariju\'an is
\begin{equation}\label{eq:divisor-complex-formula}
\operatorname{Tor}_p^P(R,\Bbbk)_e
 \cong\widetilde H_{p-1}(\Delta_e;\Bbbk)
\end{equation}
\cite[Lemma~1.1 and Theorem~1.2]{CM91}.  Their polynomial ring \(A\),
semigroup algebra \(R\), degree \(m\), and generators \(b_0,\ldots,b_g\)
are our \(P,R,e\), and \(n_1,\ldots,n_4\), with their \(g=3\).
The hypotheses are a numerical semigroup with its minimal generators
and an arbitrary coefficient field, all present here.  The indexing
in \eqref{eq:divisor-complex-formula} can also be read directly from
the degree-\(e\) part of \(R\otimes_P K_\bullet(x_1,\ldots,x_4)\):
its degree-\(p\) basis is indexed by the \(p\)-element faces of
\(\Delta_e\), and its differential is the augmented simplicial
boundary.  Hence Tor degree \(p\) corresponds to reduced homology
degree \(p-1\).

Take \(e=n_l+n_i\).  The edge \(\{l,i\}\) is a face, but no edge
joins \(\{l,i\}\) to \(\{j,k\}\): for instance \(\{l,j\}\) would
require \(n_i-n_j\in S\), impossible for distinct minimal generators.
The same argument excludes the other three cross edges.  The
\(1\)-skeleton is therefore contained in the union of the two disjoint
edges \(\{l,i\}\) and \(\{j,k\}\), and is a forest.  It has no
first homology over any field, so \(\beta_{2,n_l+n_i}^P(R)=0\).
By \eqref{eq:type3-tor2-injection} and \eqref{eq:type3-duality},
\[
\beta_{1,q+\deg M}^P(B)
 =\beta_{2,n_l+n_i}^P(B)
 \le\beta_{2,n_l+n_i}^P(R)=0.
\]
\end{proof}

We can now test all three classes in \(\alpha_1\).  The preceding
lemmas handle only the class corresponding to \(M\); the other two
still require the explicit use of the canonical-ideal membership
criterion.  Together these two mechanisms will prove that the entire
comparison map vanishes.

\begin{lemma}[Preservation of the generators of \(J\)]
\label{lem:type3-comparison}
One has \(\rho=\mu_P(J)+1\).
\end{lemma}

\begin{proof}
Since \(\varphi\) has positive degree, \(\alpha_0=0\).
The bottom of the long exact Tor sequence therefore gives
\[
\operatorname{Tor}_1^P(A_Q(-q),\Bbbk)
 \xrightarrow{\alpha_1}\operatorname{Tor}_1^P(B,\Bbbk)
 \longrightarrow\operatorname{Tor}_1^P(R,\Bbbk)
 \longrightarrow\Bbbk(-q)\longrightarrow0,
\]
and thus
\begin{equation}\label{eq:type3-cancellation}
\rho=\mu_P(J)+1-\operatorname{rank}_{\Bbbk}\alpha_1.
\end{equation}
It remains to prove \(\alpha_1=0\), rather than merely the ideal
identity \(I_S=J+(\varphi)\).

To identify this map, lift multiplication by \(\varphi\) to a homogeneous
chain map from \(K_\bullet(x_l,x_i,M)(-q)\) to a minimal resolution
\(G_\bullet\) of \(B\), taking multiplication by \(\varphi\) in degree
zero.  The chain-map identity says that the degree-one basis vector
corresponding to \(a\in\{x_l,x_i,M\}\) maps to a vector whose
\(G_1\)-boundary is \(a\varphi\).  Under the identifications
\[
\operatorname{Tor}_1^P(A_Q(-q),\Bbbk)
 \cong(Q/\mathfrak nQ)(-q),
\qquad
\operatorname{Tor}_1^P(B,\Bbbk)\cong J/\mathfrak nJ,
\]
reduction modulo \(\mathfrak n\) therefore gives
\begin{equation}\label{eq:type3-alpha1-classes}
[x_l]\longmapsto[x_l\varphi],\qquad
[x_i]\longmapsto[x_i\varphi],\qquad
[M]\longmapsto[M\varphi].
\end{equation}

For \(r=l,i\), formula~\eqref{eq:type3-C-membership} gives
\(\mathsf z^{g+n_r}\in C=QR\): the complementary integer is \(F-g\)
for \(r=l\), and \(f_i-g\) for \(r=i\), both outside \(S\).
Choose a monomial \(H_r\in Q\) of degree \(g+n_r\).  Since
\(g+n_j=\lambda n_k\), the identity
\begin{equation}\label{eq:type3-binomial-nullity}
x_r\varphi
=x_k(x_rx_k^\lambda-x_jH_r)
 +x_j(x_kH_r-x_rL)
\end{equation}
has both parenthesized binomials in \(I_S\cap Q=J\): their terms
have equal weighted degree and belong to the monomial ideal \(Q\).
Thus \(x_r\varphi\in\mathfrak nJ\) for both \(r=l,i\).
The remaining class \([M\varphi]\) has degree \(q+\deg M\) in
\(J/\mathfrak nJ\), and is zero by
Lemma~\ref{lem:type3-complementary-vanishing}.
All three images in \eqref{eq:type3-alpha1-classes} vanish, so
\(\alpha_1=0\).  Equation~\eqref{eq:type3-cancellation} proves the
assertion.
\end{proof}

\subsection{\texorpdfstring{Determining \(\mu_P(J)\)}
{Determining the number of generators of J}}

Recall that \(J=I_S\cap Q\) is the defining ideal of the auxiliary
Gorenstein quotient \(B=P/J\).  By Lemma~\ref{lem:type3-comparison},
\(\rho=\mu_P(J)+1\), so the type-three case is reduced to proving that
\(\mu_P(J)=5\).  We first restrict \(\mu_P(J)\) to \(3\) or \(5\), and then
exclude the former possibility.  The Pfaffian theorem bounds this
number once we find an element of \(J\) outside \(\mathfrak n^3\).
The relevant order is ordinary total degree, not weighted degree: a
term such as \(x_ix_j\) has order two regardless of \(n_i+n_j\).
The NG condition supplies exactly such a term in the next lemma.

\begin{lemma}[A squarefree quadratic relation in \(J\)]
\label{lem:type3-pfaffian}
The ideal \(J\) contains a binomial with a squarefree quadratic term.
Consequently \(\mu_P(J)\in\{3,5\}\).
\end{lemma}
\begin{proof}
Suppose first that \(M\ne x_j\).  Since
\(\deg M=\theta n_j+\nu n_k=n_l+F-g\),
\[
n_j+g-F=n_l-(\theta-1)n_j-\nu n_k\notin S.
\]
Indeed, since \(M\ne x_j\), the integer on the right is strictly smaller
than \(n_l\).  If it belonged to \(S\), every generator occurring in one
of its factorizations would be smaller than \(n_l\), hence would have index
in \(\{j,k\}\), because \(i>l\).  Adding
\((\theta-1)n_j+\nu n_k\) would then express \(n_l\) in
\(\langle n_j,n_k\rangle\), contrary to the minimality of \(n_l\).  If \(f_j=g\),
then \eqref{eq:NG-condition} with \(r=j\) and \(p=F\) would give
\(n_j+g-F\in S\).  Hence \(f_j=F\) or \(f_i\).  In the first case,
\eqref{eq:NG-condition} with \(p=f_i\) gives \(n_j+d\in S\); in the
second, the same condition with \(p=F\) gives \(n_j-d\in S\).
Choosing a monomial \(U\) of the corresponding degree yields
\begin{equation}\label{eq:type3-quadratic-Mneqxj}
x_ix_j-x_lU\in J
\quad\text{or}\quad
x_lx_j-x_iU\in J.
\end{equation}
Indeed, in either case the two terms have the same \(S\)-degree and
belong to \(Q\), hence the binomial lies in \(I_S\cap Q=J\).  Its
quadratic term is squarefree.

Now suppose that \(M=x_j\).  Put \(E=\{l,i,j\}\) and
\(\gamma=F+n_l\).  Since \(n_j=n_l+F-g\) and
\(f_i=F-n_i+n_l\),
\(\operatorname{PF}(S)=\{\gamma-n_r:r\in E\}\).
Write \(f_k=\gamma-n_s\) for some \(s\in E\).  For
\(t\in E\setminus\{s\}\), applying \eqref{eq:NG-condition} with
\(r=k\) and \(p=\gamma-n_t\) gives
\(n_k+n_t-n_s\in S\).  If \(H_t\) is a monomial of this degree, then
\begin{equation}\label{eq:type3-quadratic-Meqxj}
x_kx_t-x_sH_t\in J,
\end{equation}
because the two terms have the same \(S\)-degree and both lie in
\(Q=(x_l,x_i,x_j)\).  Since \(s\ne t,k\), the term \(x_kx_t\) is
squarefree and does not cancel.  Thus in either case
\(J\not\subseteq\mathfrak n^3\).

By Proposition~\ref{prop:type3-gorenstein},
\(J_{\mathfrak n}\) is a grade-three Gorenstein ideal of the regular
local ring \(P_{\mathfrak n}\).  Moreover,
\(\mu_{P_{\mathfrak n}}(J_{\mathfrak n})=\mu_P(J)\), and hence
\(\mu_P(J)\ge3\).  By the Buchsbaum--Eisenbud structure theorem,
\(\mu_P(J)\) is odd \cite[Corollary~2.2]{BE77}.  If
\(\mu_P(J)=2m+1\ge7\), then \(m\ge3\), and
\cite[Theorem~2.1]{BE77} realizes \(J_{\mathfrak n}\) as the ideal
of submaximal Pfaffians of a \((2m+1)\)-by-\((2m+1)\) alternating
matrix with entries in \(\mathfrak nP_{\mathfrak n}\).
The theorem is stated for a Noetherian local ring and a grade-three
Gorenstein ideal; these hypotheses hold here because \(P_{\mathfrak n}\)
is regular and its Cohen--Macaulay Gorenstein quotient has dimension
one.  Its minimality condition is precisely that the matrix entries
lie in the maximal ideal.  No restriction on the field characteristic
is imposed: alternating includes a zero diagonal, also in
characteristic two.  Hence
\(J_{\mathfrak n}\subseteq\mathfrak n^3P_{\mathfrak n}\), since each
such Pfaffian is a sum of products of \(m\ge3\) matrix entries.
As \(\mathfrak n^3\) is \(\mathfrak n\)-primary, contraction gives
\(J\subseteq\mathfrak n^3\), a contradiction.  Therefore
\(\mu_P(J)<7\), and consequently
\(\mu_P(J)\in\{3,5\}\).
\end{proof}

To treat the remaining case \(\mu_P(J)=3\), we first record that
\(J\) has a minimal homogeneous binomial generating set.  For each
degree \(e\), the monomials in \(Q_e\) form a \(\Bbbk\)-basis and all
map to \(\mathsf z^e\) in \(R\).  Hence
\(J_e=\ker(Q_e\to R_e)\) is spanned by binomial differences.
Therefore \(J/\mathfrak nJ\) has a homogeneous basis represented by
binomials, and graded Nakayama gives the assertion.

The Pfaffian theorem does not exclude three generators.  To test this
case, localize at a minimal prime \(\mathfrak p\) of \(Q\).
The equation \(\varphi\) becomes a unit, so \(J_{\mathfrak p}=Q_{\mathfrak p}\).
Three global generators of \(J\) would therefore have to form a basis
of the three-dimensional vector space
\[
Q_{\mathfrak p}/\mathfrak pQ_{\mathfrak p}
 \cong(Q/Q^2)\otimes_{P/Q}\kappa(\mathfrak p).
\]
In particular, none could have zero image.  If \(M\) involves both
remaining variables, testing this condition at both minimal primes
will be incompatible with the minimal generators of \(C\).  If \(M\)
is a pure power, the same basis condition forces a permutation of the
three monomial generators of \(Q\).  These are the two steps of the
next lemma.

\begin{lemma}[Structure when \(\mu_P(J)=3\)]
\label{lem:type3-ci-form}
Suppose \(\mu_P(J)=3\).  Then \(\nu=0\).  Put
\(E=\{l,i,j\}\) and
\(\varepsilon_l=\varepsilon_i=1\), \(\varepsilon_j=\theta\).
There exist a permutation \(\sigma:E\to E\) and integers
\(A_r>\varepsilon_r\), \(u_r\ge1\) such that
\begin{equation}\label{eq:type3-ci-form}
J=\left(
 x_r^{A_r}-x_{\sigma(r)}^{\varepsilon_{\sigma(r)}}x_k^{u_r}
 :r\in E
\right).
\end{equation}
\end{lemma}

\begin{proof}
Choose three homogeneous binomial generators of \(J\).  Let
\(\mathfrak p\) be a minimal prime of \(Q\).  Since
\(\bar\varphi\) is a nonzerodivisor on \(P/Q\), one has
\(\varphi\notin\mathfrak p\).  As \(\varphi\in I_S\), it becomes a
unit in \(P_{\mathfrak p}\), and hence
\((I_S)_{\mathfrak p}=P_{\mathfrak p}\) and
\(J_{\mathfrak p}=Q_{\mathfrak p}\).
Since \(Q\) is a height-three complete intersection,
\(Q_{\mathfrak p}\) has three minimal generators.  Thus the localized
classes of our three binomial generators form a basis of
\(Q_{\mathfrak p}/\mathfrak pQ_{\mathfrak p}\) over
\(\kappa(\mathfrak p)\); in particular, none of them vanishes.

\emph{Testing both minimal primes forces \(\nu=0\).}
Suppose \(\nu>0\).  The two minimal
primes of \(Q\) are
\(\mathfrak p_j=(x_l,x_i,x_j)\) and
\(\mathfrak p_k=(x_l,x_i,x_k)\).
Specializing \(x_i=x_j=x_k=0\) does not annihilate \(J\), since
\(x_l(x_l^{n_k}-x_k^{n_l})\in J\).  Hence one of the three generators
has the form \(x_l^A-N\), where \(A\ge2\) by the minimality of \(n_l\).
The class of \(x_l^A\) vanishes in both
\(Q_{\mathfrak p_j}/\mathfrak p_jQ_{\mathfrak p_j}\) and
\(Q_{\mathfrak p_k}/\mathfrak p_kQ_{\mathfrak p_k}\); therefore the
class of \(N\) is nonzero in both quotients.

Since \(x_k\) is a unit in \(P_{\mathfrak p_j}\), the classes of
\(x_l,x_i,x_j^\theta\) form a basis of
\(Q_{\mathfrak p_j}/\mathfrak p_jQ_{\mathfrak p_j}\).  Thus a monomial
has nonzero class there precisely when it is
\(x_lx_k^b\), \(x_ix_k^b\), or \(x_j^\theta x_k^b\), with \(b\ge0\).
Similarly, since \(x_j\) is a unit in \(P_{\mathfrak p_k}\), the
corresponding monomials are
\(x_lx_j^a\), \(x_ix_j^a\), and \(x_k^\nu x_j^a\), with \(a\ge0\).
The common possibilities are therefore \(N=x_l,x_i,M\).
The first gives \(A=1\) by homogeneity, while either of the other two
makes \(x_i\), respectively \(M\), redundant in the minimal generating
set \(C=(x_l,x_i,M)R\).  All three are impossible, so \(\nu=0\).

\emph{Testing the remaining prime determines the generators.}
Thus
\(Q=(x_l,x_i,x_j^\theta)\) has the unique minimal prime
\(\mathfrak p=(x_l,x_i,x_j)\).  For each \(r\in E\), specialization
at \(x_s=0\) for \(s\ne r\) does not annihilate \(J\), as witnessed by
\(x_r^{\varepsilon_r}(x_r^{n_k}-x_k^{n_r})\in J\).
Hence some generator has a pure-power term \(x_r^{A_r}\).
Since \(Q\) is monomial, \(A_r\ge\varepsilon_r\).  If
\(A_r=\varepsilon_r\), write the corresponding binomial as
\(x_r^{\varepsilon_r}-N\).  Then \(N\in Q\), so
\(x_s^{\varepsilon_s}\mid N\) for some \(s\in E\).  The case \(s=r\)
is impossible by equality of weighted degrees, while \(s\ne r\)
would express the generator \(x_r^{\varepsilon_r}\) of \(C\) as a
multiple of the distinct generator \(x_s^{\varepsilon_s}\).
Thus \(A_r>\varepsilon_r\).

Consequently every pure-power term \(x_r^{A_r}\) vanishes in
\(Q_{\mathfrak p}/\mathfrak pQ_{\mathfrak p}\).  No binomial generator
can contain two such terms, since its class would then vanish.  There
are three required pure-power terms and three generators, so the
generators may be indexed by \(r\in E\) according to their pure-power
terms \(x_r^{A_r}\).  The other term of each generator must have
nonzero class in \(Q_{\mathfrak p}/\mathfrak pQ_{\mathfrak p}\), and
hence has the form
\(x_s^{\varepsilon_s}x_k^u\), with \(s\in E\) and \(u\ge0\).
If \(u=0\), then \(s=r\) contradicts
\(A_r>\varepsilon_r\), while \(s\ne r\) contradicts the minimality of
the generators of \(C\).  Hence \(u\ge1\).

Finally, \(x_k\) is a unit in \(P_{\mathfrak p}\), so the class of
\(x_s^{\varepsilon_s}x_k^u\) is a nonzero scalar multiple of the basis
class \(x_s^{\varepsilon_s}\).  Since the three localized generator
classes are linearly independent, their corresponding indices \(s\)
are distinct.  They therefore define a permutation
\(\sigma:E\to E\), and
\[
J=\left(
x_r^{A_r}
-x_{\sigma(r)}^{\varepsilon_{\sigma(r)}}x_k^{u_r}
:r\in E
\right),
\]
with \(A_r>\varepsilon_r\) and \(u_r\ge1\), as required.
\end{proof}

The localized basis condition has now forced a global binomial form.
To contradict it we need an ideal-membership test, not just a comparison
of numerical degrees: distinct monomials can have the same weighted
degree without being connected by the three proposed relations.  The
following elementary criterion records the precise test used in the
last case of the proof.

\begin{lemma}[Degree-wise binomial graph]
\label{lem:binomial-graph}
Let
\(G_r=U_r-V_r\), \(1\le r\le m\), be homogeneous binomials and put
\(L=(G_1,\ldots,G_m)\).  For a fixed degree \(e\), let
\(\Gamma_e\) be the graph whose vertices are the monomials of degree
\(e\), with an edge between \(WU_r\) and \(WV_r\) whenever \(W\) is
a monomial and these products have degree \(e\).  If
\(U-V\in L_e\) for monomials \(U,V\) of degree \(e\), then \(U\)
and \(V\) lie in the same connected component of \(\Gamma_e\).
\end{lemma}

\begin{proof}
The vector space \(L_e\) is spanned by the edge differences
\(WU_r-WV_r\).  For each connected component \(\mathcal C\) of
\(\Gamma_e\), define a linear functional on the degree-\(e\) part of
\(P\) by summing the coefficients of the monomials belonging to
\(\mathcal C\).  This functional vanishes on every edge difference and
therefore on \(L_e\).  If \(U\) and \(V\) belonged to different
components, applying the functional for the component of \(U\) to
\(U-V\) would give \(1\), contrary to \(U-V\in L_e\).
\end{proof}

\begin{proposition}\label{prop:type3}
If \(\tau=3\), then \(\rho=6\).
\end{proposition}

\begin{proof}
By Lemmas~\ref{lem:type3-comparison} and~\ref{lem:type3-pfaffian},
it suffices to exclude \(\mu_P(J)=3\).
Suppose this equality holds, and use
Lemma~\ref{lem:type3-ci-form}.  In particular \(M=x_j^\theta\).

If \(\theta>1\), reducing \eqref{eq:type3-ci-form} modulo \(x_k\)
gives the monomial ideal
\begin{equation*}
\frac{J+(x_k)}{(x_k)}
=(x_l^{A_l},x_i^{A_i},x_j^{A_j})
\subseteq\Bbbk[x_l,x_i,x_j].
\end{equation*}
The binomial in \eqref{eq:type3-quadratic-Mneqxj} has the
squarefree term \(x_ix_j\) or \(x_lx_j\), which survives this
reduction and is divisible by none of these pure powers.
Its other term cannot cancel it, because it contains the other of
\(x_l,x_i\).  This contradicts membership in the displayed monomial
ideal.

It remains that \(\theta=1\).  We now use the NG condition to obtain
two binomials involving a fixed variable \(x_s\).  The normal form
\eqref{eq:type3-ci-form} will force both to use the same generator of
\(J\), but with different other endpoints in the degree-wise graph.
Choose \(s\in E\) as in
\eqref{eq:type3-quadratic-Meqxj}.  For both distinct indices
\(t\in E\setminus\{s\}\), one has
\begin{equation*}
x_kx_t-x_sH_t\in J.
\end{equation*}
Let \(e_t=\deg(x_kx_t)\), and form the graph \(\Gamma_{e_t}\)
from Lemma~\ref{lem:binomial-graph} using the three generators in
\eqref{eq:type3-ci-form}.  Since
\(x_kx_t-x_sH_t\in J_{e_t}\), that lemma shows that
\(x_kx_t\) and \(x_sH_t\) lie in the same connected component.

No pure power \(x_r^{A_r}\), with \(A_r\ge2\), divides \(x_kx_t\).
The only generator in \eqref{eq:type3-ci-form} whose nonpure term can
divide \(x_kx_t\) is therefore the unique one indexed by \(r\) with
\(\sigma(r)=t\); such divisibility is possible only when \(u_r=1\),
and then the corresponding edge is
\[
x_kx_t\;\longleftrightarrow\;x_r^{A_r}.
\]
This edge must occur, because \(x_sH_t\ne x_kx_t\) and the two
monomials lie in the same component.  Starting from
\(x_r^{A_r}\), no generator other than the one indexed by \(r\) can
produce an edge: no other pure-power term divides \(x_r^{A_r}\), and
every nonpure term in \eqref{eq:type3-ci-form} contains \(x_k\).
For the generator indexed by \(r\), the only possible multiplier is
\(1\).  Hence the connected component of \(x_kx_t\) consists exactly
of
\[
\{x_kx_t,x_r^{A_r}\}.
\]
It follows that \(x_sH_t=x_r^{A_r}\).  Equality of these monomials
forces \(r=s\), while the edge above then gives \(\sigma(s)=t\).
There are two distinct choices of \(t\in E\setminus\{s\}\), which is
impossible for the function \(\sigma\).

Consequently \(\mu_P(J)=5\), and
Lemma~\ref{lem:type3-comparison} gives \(\rho=6\).
\end{proof}

\begin{proof}[Proof of Theorem~\ref{thm:main}]
By Section~\ref{sec:preliminaries}, one has \(\tau\in\{2,3\}\).
If \(\tau=2\), Proposition~\ref{prop:type2} gives
\(\rho\in\{4,5\}\); if \(\tau=3\), Proposition~\ref{prop:type3}
gives \(\rho=6\).  Substituting these possibilities into the Betti
formula \((1,\rho,\rho+\tau-1,\tau)\) from
Proposition~\ref{prop:betti-identities} gives the three asserted
Betti sequences.  Proposition~\ref{prop:sharpness} shows that each of
them occurs.
\end{proof}

\appendix
\section{Sharpness examples}\label{app:sharpness}

We verify that all three Betti sequences in Theorem~\ref{thm:main} are
realized.  The calculation is independent of the coefficient field.  For
\(b\in S\), let \(\mathcal G_b\) be the graph whose vertices are the
factorizations
\[
(a_1,a_2,a_3,a_4)\in\mathbb N^4,
\quad \sum_{r=1}^4a_rn_r=b,
\]
and in which two vertices are adjacent when their supports intersect.  Here adjacency is defined by shared support of arbitrary factorizations,
not by the chosen binomial moves in Lemma~\ref{lem:binomial-graph}.
The connected components are the usual \(R\)-classes of the
factorization fibre.  Let \(c_b\) denote their number, and define
\[
\operatorname{Betti}(S)=\{b\in S\setminus\{0\}:c_b>1\}.
\]
The factorization-graph description of minimal presentations gives
\[
\rho=\sum_{b\in\operatorname{Betti}(S)}(c_b-1).
\]
Indeed, Rosales works with a subsemigroup of \(\mathbb N\), a fixed
generating set, the factorization map \(\phi\), and the equivalence
classes generated by nonzero scalar products of factorizations
\cite[Theorem~1.2 and Corollary~1.3]{Rosales96}.
In the present notation his fibre \([b]=\phi^{-1}(b)\) is our vertex
set, and his classes are exactly the components of \(\mathcal G_b\).
His criterion says that a minimal presentation is obtained by choosing
a tree on these components for each disconnected fibre, using exactly
\(c_b-1\) relations.  Our positive minimal generators make every fibre
finite.  Herzog's correspondence then identifies those relations with
minimal binomial generators over any field
\cite[Proposition~1.5]{Herzog70}.

\begin{proposition}[Sharpness]\label{prop:sharpness}
Each of the three Betti sequences in Theorem~\ref{thm:main} occurs over every
field.  More precisely,
\[
\begin{array}{c|c|c|c}
S & \operatorname{PF}(S) & \text{an NG-vector} & \rho \\ \hline
\langle6,7,8,9\rangle & \{10,11\} & (11,10,10,10) & 4\\
\langle6,11,13,14\rangle & \{16,21\} & (21,16,21,21) & 5\\
\langle5,6,9,13\rangle & \{4,7,8\} & (8,7,4,4) & 6.
\end{array}
\]
\end{proposition}

\begin{proof}
Write the three semigroups as \(S_4,S_5,S_6\), according to the value of
\(\rho\) in the table.  Their gaps up to the Frobenius number are respectively
\[
\begin{aligned}
\mathbb N\setminus S_4&=\{1,2,3,4,5,10,11\},\\
\mathbb N\setminus S_5&=\{1,2,3,4,5,7,8,9,10,15,16,21\},\\
\mathbb N\setminus S_6&=\{1,2,3,4,7,8\}.
\end{aligned}
\]
The consecutive blocks \(12,\ldots,17\), \(22,\ldots,27\), and
\(9,\ldots,13\) belong to \(S_4,S_5,S_6\), respectively.
Their lengths are the corresponding multiplicities \(6,6,5\); adding
that multiplicity proves that every integer beyond each block also
belongs to the semigroup.  The displayed gap lists are therefore
exhaustive, and the pseudo-Frobenius condition can be checked on them.  It gives exactly the three pseudo-Frobenius sets shown
in the statement.  Substituting those sets into \(n_r+f_r-p\in S\) verifies
the displayed NG-vectors coordinate by coordinate, so all three semigroups
are nearly Gorenstein.  None is almost symmetric: the constant Frobenius
vector fails respectively at
\[
9+11-10=10\notin S_4,
\quad
11+21-16=16\notin S_5,
\quad
6+8-7=7\notin S_6.
\]

It remains to compute the presentation cardinalities.  We first display the
disconnected fibres used below; braces on the right separate the connected
components of the corresponding graph.  We then give a terminating reduction argument proving that the resulting binomials generate the full toric
ideal, and hence that no disconnected fibre has been omitted.  For
\(S_4=\langle6,7,8,9\rangle\),
\[
\begin{array}{c|c}
14 & \{(0,2,0,0)\}\;\sqcup\;\{(1,0,1,0)\}\phantom{.}\\
15 & \{(0,1,1,0)\}\;\sqcup\;\{(1,0,0,1)\}\phantom{.}\\
16 & \{(0,0,2,0)\}\;\sqcup\;\{(0,1,0,1)\}\phantom{.}\\
18 & \{(0,0,0,2)\}\;\sqcup\;\{(3,0,0,0)\}.
\end{array}
\]
The four relations joining the components supplied by these fibres are
\[
x_2^2-x_1x_3,\quad
x_2x_3-x_1x_4,\quad
x_3^2-x_2x_4,\quad
x_4^2-x_1^3.
\]
For \(S_5=\langle6,11,13,14\rangle\), the displayed disconnected fibres are
\[
\begin{array}{c|c}
24 & \{(0,1,1,0)\}\;\sqcup\;\{(4,0,0,0)\}\phantom{.}\\
25 & \{(0,1,0,1)\}\;\sqcup\;\{(2,0,1,0)\}\phantom{.}\\
26 & \{(0,0,2,0)\}\;\sqcup\;\{(2,0,0,1)\}\phantom{.}\\
28 & \{(0,0,0,2)\}\;\sqcup\;\{(1,2,0,0)\}\phantom{.}\\
33 & \{(0,3,0,0)\}\;\sqcup\;\{(1,0,1,1)\}.
\end{array}
\]
The corresponding five relations joining the components are
\[
\begin{gathered}
x_2x_3-x_1^4,\quad x_2x_4-x_1^2x_3,\quad x_3^2-x_1^2x_4,\\
x_4^2-x_1x_2^2,\quad x_2^3-x_1x_3x_4.
\end{gathered}
\]
Finally, for \(S_6=\langle5,6,9,13\rangle\), the displayed disconnected
fibres are
\[
\begin{array}{c|c}
15 & \{(0,1,1,0)\}\;\sqcup\;\{(3,0,0,0)\}\\
18 & \{(0,0,2,0)\}\;\sqcup\;\{(0,3,0,0)\}\;\sqcup\;\{(1,0,0,1)\}\\
19 & \{(0,1,0,1)\}\;\sqcup\;\{(2,0,1,0)\}\\
22 & \{(0,0,1,1)\}\;\sqcup\;\{(2,2,0,0)\}\\
26 & \{(0,0,0,2)\}\;\sqcup\;\{(1,2,1,0),(4,1,0,0)\}.
\end{array}
\]
The degree \(18\) fibre has three components and therefore contributes two
minimal relations; each of the other four Betti fibres contributes one.  A
minimal binomial generating set is, for example,
\begin{equation}\label{eq:sharpness-type3-generators}
\begin{aligned}
h_1&=x_2x_3-x_1^3,       & h_2&=x_3^2-x_2^3,\\
h_3&=x_2^3-x_1x_4,       & h_4&=x_2x_4-x_1^2x_3,\\
h_5&=x_3x_4-x_1^2x_2^2, & h_6&=x_4^2-x_1x_2^2x_3.
\end{aligned}
\end{equation}
We now certify that the displayed relations generate the entire defining
ideal in each case.  Let \(G_r\) be the displayed set of binomials for
\(S_r\), and put \(L_r=(G_r)\subseteq P\).  Every element of \(G_r\)
has equal \(S_r\)-degree on its two terms, so \(L_r\subseteq I_{S_r}\).
Orient every displayed relation from its left-hand monomial to its
right-hand monomial.  These rewrites strictly decrease the auxiliary
monomial weights
\[
\mathbf w_4=(1,3,3,2),\quad
\mathbf w_5=(1,4,4,5),\quad
\mathbf w_6=(1,3,5,7),
\]
where the coordinates correspond to \((x_1,x_2,x_3,x_4)\).  Hence every
reduction terminates.  A monomial is irreducible for the respective systems
precisely when it has the form \(x_1^qu\), with \(q\ge0\) and
\[
\begin{aligned}
\mathcal B_4&=\{1,x_2,x_3,x_4,x_2x_4,x_3x_4\},\\
\mathcal B_5&=\{1,x_2,x_2^2,x_3,x_4,x_3x_4\},\\
\mathcal B_6&=\{1,x_2,x_2^2,x_3,x_4\}.
\end{aligned}
\]
For \(m\in S\setminus\{0\}\), write
\(\operatorname{Ap}(S,m)=\{s\in S:s-m\notin S\}\) for the Ap\'ery set.
The \(S_r\)-degrees of the residual monomials are
\[
\begin{aligned}
\{\deg u:u\in\mathcal B_4\}
 &=\{0,7,8,9,16,17\}=\operatorname{Ap}(S_4,6),\\
\{\deg u:u\in\mathcal B_5\}
 &=\{0,11,13,14,22,27\}=\operatorname{Ap}(S_5,6),\\
\{\deg u:u\in\mathcal B_6\}
 &=\{0,6,9,12,13\}=\operatorname{Ap}(S_6,5).
\end{aligned}
\]
Indeed, in each line the displayed degrees represent all residue classes
modulo the first generator, and subtracting that generator from every
nonzero degree gives one of the gaps listed above.  Thus every element of
\(S_r\) has a unique expression
\[
qn_1+\deg u,
\quad q\ge0,\quad u\in\mathcal B_r.
\]
The induced surjection
\[
P/L_r\longrightarrow\Bbbk[S_r]
\]
therefore sends a spanning set of irreducible monomials bijectively to the
monomial basis of \(\Bbbk[S_r]\).  It is an isomorphism, and consequently
\(L_r=I_{S_r}\).

The chosen binomials join distinct components in every displayed fibre and,
for \(S_6\) in degree \(18\), form a tree on its three components.  Thus no
generator is redundant, so \(G_r\) is a minimal binomial generating set,
equivalently a minimal presentation of \(S_r\).  Rosales's criterion then
excludes any additional disconnected fibre, since it would have to
contribute a relation in a further degree.  Accordingly
\[
\rho(S_4)=4,\quad \rho(S_5)=5,\quad \rho(S_6)=6.
\]
Combining these values with the types \(2,2,3\) and
Proposition~\ref{prop:betti-identities} gives respectively
\[
(1,4,5,2),\quad (1,5,6,2),\quad (1,6,8,3).
\]
This proves the proposition.
\end{proof}

\section*{Declaration of AI use}

OpenAI's ChatGPT was used as an assistive tool during the development of this
work, including for literature searches, exploration and checking of mathematical
arguments, and manuscript preparation. The author led the overall research
direction, determined the conceptual framework and proof architecture, and made
all substantive mathematical decisions. All AI-generated suggestions were
independently evaluated by the author. The author independently checked all
mathematical claims, proofs, and citations, determined the final form of the
manuscript, and takes full responsibility for the content of the work.

\end{document}